\documentclass[a4paper,10pt]{article}
\usepackage{amssymb}
\usepackage{amsmath,amsfonts,amsthm,amssymb}
\usepackage[dvips]{graphics}
\usepackage{epsfig}
\usepackage{indentfirst}
\usepackage{color}
\allowdisplaybreaks

\newtheoremstyle{theorem}
  {10pt}          
  {10pt}  
  {\sl}  
 {}
  {\bf}  
  {. }    
  { }    
  {}     
\theoremstyle{theorem}

\newtheorem{theorem}{Theorem}[section]

\newtheorem{remark}{Remark}[section]

\numberwithin{equation}{section}

\newtheoremstyle{defi}
{10pt}  
{10pt}  
{\rm}   
{}      
{\bf}   
{. }    
{ }     
{}      
\theoremstyle{defi}

\begin{document}
\baselineskip = 13.5pt

\title{\bf Low Mach and thin domain limit for the compressible Euler system}

\author{Matteo Caggio$^{1}$ \ \ \ Bernard Ducomet$^{2}$ \ \ \ \v S\' arka Ne\v casov\' a$^{1}$ \ \ \  Tong Tang$^{3}$
\\
{\small  1. Institute of Mathematics of the Academy of Sciences of the Czech Republic,} \\
{\small \v Zitn\' a 25, 11567, Praha 1, Czech Republic}\\
{\small 2. Universit\'e Paris-Est, LAMA, 61 Avenue du G\'en\'eral de Gaulle, F-94010 Cr\'eteil, France}\\
{\small 3. Department of Mathematics, College of Sciences,}\\
{\small Hohai University, Nanjing 210098, P.R. China}\\
\date{}}
\maketitle
\begin{abstract}
We consider the compressible Euler system describing the motion of an ideal fluid confined to a straight layer $\Omega_{\delta}=(0,\delta)\times\mathbb{R}^2, \ \ \delta>0$. In the framework of \textit{dissipative measure-valued solutions}, we show the convergence to the strong solution of the 2D incompressible Euler system when the Mach number tends to zero and $\delta\rightarrow0$.
\vspace{0.5cm}

{{\bf Key words:} compressible Euler equations, dissipative measure-valued solutions, low Mach number, thin domain.}

\medskip

\end{abstract}
\maketitle
\section{Introduction}\setcounter{equation}{0}
The present paper is devoted to the problem of the limit passage from three-dimensional to two-dimensional geometry, and from compressible to incompressible inviscid fluid. In the infinite straight layer
\begin{equation} \label{Omega}
\Omega_{\delta}=(0,\delta)\times\mathbb{R}^2, \ \ \delta>0,
\end{equation}
we consider the compressible Euler system describing the motion of a barotropic fluid,
\begin{equation} \label{cont}
\partial_{t}\varrho_{\epsilon}+\textrm{div}_{x}\left(\varrho_{\epsilon}\mathbf{u}_{\epsilon}\right)=0,
\end{equation}
\begin{equation} \label{mom}
\partial_{t}\left(\varrho_{\epsilon}\mathbf{u}_{\epsilon}\right)+\textrm{div}_{x}\left(\varrho_{\epsilon}\mathbf{u}_{\epsilon}\otimes\mathbf{u}_{\epsilon}\right)+\frac{1}{\epsilon^{2}}\nabla_{x}p\left(\varrho_{\epsilon}\right)=0,
\end{equation}
supplemented with the initial conditions
\begin{equation} \label{ic}
\rho_\epsilon (0, \cdot) = \rho_{0, \epsilon}, \ \
\rho_\epsilon\mathbf{u}_\epsilon (0, \cdot) = \mathbf{m}_{0, \epsilon},
\end{equation}
and the far field conditions
\begin{equation} \label{ac}
\mathbf{u}_{\varepsilon}\rightarrow 0, \ \ \rho_{\epsilon}\rightarrow{\widetilde{\rho}} \ \ \textrm{as} \ \ \left|x\right|\rightarrow\infty,
\end{equation}
with $\widetilde{\rho}>0$ constant.

The above system is written in its non-dimensional form, with $\epsilon$ the Mach number (the Strouhal number allowed by the scale analysis is set equal to one). Here, $\varrho_{\epsilon}=\varrho_{\epsilon}\left(x,t\right)$  $\mathbf{u}_{\epsilon}=\mathbf{u}_{\epsilon}\left(x,t\right)$ and $p=p(\varrho_{\epsilon}\left(x,t\right))$ represents the mass density, the velocity vector and the pressure respectively (for assumptions on the pressure see Theorem \ref{TH}).

In the context of the low Mach number limit, the convergence of the solution of the compressible Euler system to the solution of the incompressible system was shown in several papers for \textit{well-prepared} initial data of the compressible system, namely data for which the acoustic waves are not allowed, and for smooth solutions of the compressible flow (see \cite{As, Eb, KlMa, MeSc, Uk}). Indeed, it is known that solutions of the compressible Euler system develop singularities in a finite time independently how smooth and/or small the initial data are. Moreover, it was shown by Feireisl et al. \cite{e6}, it is very hard to prove that the life span of the smooth solutions is independent by the Mach number. Consequently, smooth solutions are quite restrictive for compressible inviscid flows. We would like to stress that, under certain hypotheses on initial data, Serre and Grassin \cite{gr1,gr2} proved global smooth solutions to the compressible Euler in $\mathbb{R}^d \  (d\geq1)$. Recently, these results have been extended by Blanc et al. \cite{bl} for the Euler system coupled to the Helmholtz or Poisson equations.

In order to obtain global existence results, it is necessary to move into weak solutions. As mentioned by Feireisl et al. \cite{e6}, the recent theory of \textit{convex integration}  has shown the existence of a "large number" of global-in-time weak solutions for regular initial data, even if "most of them" violate the basic energy inequality associated to the system (see \cite{c, le}). Moreover, there exist a family of "wild" initial data producing infinitely many weak solutions satisfying the so-called \textit{admissibility criteria} (see, for example, \cite{c1, c2, le}). Despite the present results the existence of global-in-time \textit{admissible} weak solutions for arbitrary (smooth) initial data remains an open problem for the compressible Euler system. However, a recent analysis by Breit et al. \cite{brf} proposed a new concept of dissipative solutions to the compressible Euler system.

Therefore, for the purpose of the present analysis, we move into the so-called \textit{dissipative measure-valued (DMV) solutions}.
Let us mention that such type of solutions were recently developed by
Demoulini et al. \cite{De} in the context of polyconvex
elastodynamics, while in context of fluid dynamics were studied by
Gwiadza et al. \cite{e3, g} (measure-valued solutions in the sense of
DiPerna \cite{Di} can be found e.g. \cite{m}, \cite{Ne}, \cite{BBN}).
Indeed, for the Euler system, the advantage of DMV solutions is that they exist globally in time for any finite energy initial data even if, in general, the solutions are not uniquely determined by the initial data. In this framework, very recently, the low Mach number limit analysis for the compressible Euler system was developed by Feireisl et al. \cite{e6} for the case of \textit{well-prepared} and \textit{ill-prepared} initial data, namely the acoustic waves are present and they are dispersed in the low Mach number limit. For other recent results concerning singular limit analysis in the framework of DMV solutions the reader can refer, for example, to \cite{br, b}. While, for heat conducting inviscid fluid to \cite{j1, j2}.

In the context of thin-domain limit analysis, recently Caggio et al. \cite{cg} consider the compressible Navier-Stokes system describing the motion of a viscous fluid confined in the straight layer (\ref{Omega}). They show that the weak solutions in the 3D domain converge strongly to the solution of the 2D incompressible Navier-Stokes equations (Euler equations) when the Mach number tends to zero as well as $\delta\rightarrow0$ (and the viscosity goes to zero). For other works related to a dimension reduction limit the reader can refer, for example, to \cite{be, j3, du, ma}. Compared with Naver-Stokes equations, it is interesting and natural to investigate the low Mach number limit for compressible Euler system in the context of thin-domain limit analysis. Based on the analysis of Caggio et al. \cite{cg} and Feireisl et al. \cite{e6}, our goal is to rigorously justify the asymptotic limit of solutions $\left(\varrho_{\epsilon},\mathbf{u}_{\epsilon}\right)$ for $\epsilon,\delta\rightarrow0$ in the framework of DMV solutions. The analysis will be developed for ill-prepared data. The well-prepared case is an easy consequence.

Formally, for $\epsilon,\delta\rightarrow0$, namely the low Mach and thin-domain limit, the limit system is expected to satisfy the two-dimensional incompressible Euler system
\begin{equation} \label{cont_E}
\text{div}_h\mathbf{v}=0,
\end{equation}
\begin{equation} \label{mom_E}
\partial_t\mathbf{v}+\mathbf{v}\nabla_h\mathbf{v}+\nabla_h\Pi=0,
\end{equation}
with $\mathbf{v}=\mathbf{v}\left(x,t\right)$ the velocity, $\Pi=\Pi\left(x,t\right)$ the pressure, and the initial condition
\begin{equation} \label{ic_E}
\mathbf{v} (0, \cdot)=\mathbf{v}_0.
\end{equation}
Here, we meant $\mathbf{v}=(v_1,v_2)$ a vector field in $\mathbb{R}^{2}$, with
$$
\nabla_h = (\partial_{x_1}, \partial_{x_2}), \ \ \text{div}_h = \nabla_h \cdot.
$$
We assume the thickness $\delta$ of the domain depending by $\epsilon$ such that $\delta(\epsilon)\rightarrow0$ as $\epsilon$, and, for a function $f$ defined in $\Omega_{\delta}$, we denote its average in the $x_3$-direction as
\begin{equation} \label{avg}
\overline{f}(x_h) =\overline{f}^\delta (x_h) = \frac{1}{\delta}\int_0^\delta f(x_h,x_3)dx_3,
\end{equation}
where $x_h = (x_1,  x_2)$. Namely, a bar over a function denotes the average over $x_{3}\in\left(0,\delta\right)$.
We will show that, if the average value of the initial data of the compressible Euler system converges, in a certain sense specified below, to the initial data of the incompressible Euler system, namely $\left(\overline{\rho_{0,\epsilon}},\overline{\mathbf{u}_{0,\epsilon}}\right)\rightarrow\left(1,\mathbf{v}_{0}\right)$, then the limit $\left(\overline{\rho_{\epsilon}},\overline{\mathbf{u}_{\epsilon}}\right)$ -- the (space) average of the solutions $\left(\rho_{\epsilon},\mathbf{u}_{\epsilon}\right)$ of the compressible Euler system -- corresponds to the (classical) solution of the incompressible Euler system in $\mathbb{R}^{2}$.

The paper is organized as follows. In Section 2, we introduce the definition of dissipative measure solutions, relative energy and the other necessary material. In Section 3, we  perform the necessary analysis of the acoustic waves and state our main theorem. Section 4 is devoted to deriving uniform bounds of the Euler system independent of $\epsilon$. The proof of the main theorem is completed in Section 5.

\section{Measure-valued solutions and relative energy inequality}

In this section we introduce the DMV solutions and the relative energy inequality as key tool of our analysis. Smooth solution of the incompressible Euler system will be also discussed.

\subsection{Dissipative measure-valued solutions}
Before introducing the notion of DMV solutions, for convenience of readers, we give a short introduction of some concepts. For more details see \cite{j2} and \cite{CFKW}.

First, we introduce the phase-space associated to the solutions $[\rho,\mathbf{m}]=[\rho,\rho\mathbf{u}]$, namely

$$\mathcal{Q}=\{[\rho,\mathbf{m}] \ | \  \rho\in[0,\infty), \ \mathbf{m}\in\Omega_{\delta}\}.$$
Now, let $L^\infty_{weak-(*)}((0,T)\times\Omega;\mathcal{P}(\mathcal{Q}))$ be the space of essentially bounded weakly-$(*)$ measure maps $Y:(0,T)\times\Omega\rightarrow \mathcal{P}(\mathcal{Q})$, $(t,x)\mapsto Y_{t,x}$. By virtue of fundamental theorem on Young measures (see Ball \cite{Ba})  there exists a subsequence of $\{\rho_\epsilon,\mathbf m_\epsilon\}_{\epsilon>0}$ and parameterized family of probability measures $\{Y_{t,x}\}_{(t,x)\in(0,T)\times\Omega}$
\begin{eqnarray*}
[(t,x)\mapsto Y_{t,x}]\in L^\infty_{weak-(*)}((0,T)\times\Omega;\mathcal{P}(\mathcal{Q})),
\end{eqnarray*}
such that a.a. $(t,x)\in(0,T)\times\Omega$
\begin{eqnarray*}
\langle Y_{t,x}; G(\rho,\mathbf m)\rangle=\widehat{G(\rho,\mathbf m)}(t,x) \ \mbox{for any} \ G\in C_c(\mathcal{Q}), \ \mbox{and} \ \mbox{a.a.} \ (t,x)\in(0,T)\times\Omega
\end{eqnarray*}
whenever
\begin{eqnarray*}
G(\rho_\epsilon,\mathbf m_\epsilon)\rightarrow\widehat{G(\rho,\mathbf m)}(t,x) \ \mbox{weakly-$(*)$} \ \mbox{in} \ L^\infty((0,T)\times\Omega).
\end{eqnarray*}
Above, the hat over a function is intended as weak limit.
The parameterized family of measures $\{Y_{t,x}\}_{(t,x)\in(0,T)\times\Omega}$ is called \textit{Young measure} associated to the sequence
$\{\rho_\epsilon,\mathbf m_\epsilon\}_{\epsilon>0}$.
If $G\in C(\mathcal{Q})$ is such that
\begin{eqnarray*}
\int^T_0\int_\Omega |G(\rho_\epsilon,\mathbf m_\epsilon)|dx\leq C,
\end{eqnarray*}
then $G$ is $Y_{t,x}$ integrable for almost all $(t,x)\in(0,T)\times\Omega$ and
\begin{eqnarray*}
[(t,x)\mapsto \langle Y_{t,x}; G(\rho,\mathbf m)\rangle]\in L^1((0,T)\times\Omega),
\end{eqnarray*}
and
\begin{eqnarray*}
G(\rho_\epsilon,\mathbf m_\epsilon)\rightarrow\widehat{G(\rho,\mathbf m)}(t,x) \ \mbox{weakly-$(*)$ in} \ \mathcal{M}((0,T)\times\Omega).
\end{eqnarray*}
Note that the Young measure $[(t,x)\mapsto \langle Y_{t,x}; G(\rho,\mathbf m)\rangle]$ is a parameterized family of non-negative measures acting on the phase space $\mathcal{Q}$, while $\widehat{G(\rho,\mathbf m)}(t,x)$ is a signed measure on the physical space $[0,T]\times \Omega$.
In conclusion, the difference
\begin{eqnarray*}
\mu_G\equiv\widehat{G(\rho,\mathbf m)}-[(t,x)\mapsto \langle Y_{t,x}; G(\rho,\mathbf m)\rangle]\in \mathcal{M}((0,T)\times\Omega),
\end{eqnarray*}
is called \textit{concentration defect measure}.


A \emph{dissipative measure-valued (DMV) solution} of the Euler system (\ref{cont}) - (\ref{mom}) is a Young measure $\{Y_{t,x}\}_{t\in[0,T],x\in\Omega_{\delta}}$
satisfying:

\bigskip
\noindent
$\bullet$ Equation of continuity
\begin{align} \label{cont_mv}
\int^T_0\int_{\Omega_{\delta}}[\langle Y_{t,x};\rho\rangle\partial_t\varphi+\langle Y_{t,x};\mathbf{m}\rangle\nabla_x\varphi]dxdt
=-\int_{\Omega_{\delta}}\langle Y_{0,x};\rho\rangle\varphi(0)dx,
\end{align}
for all $\varphi\in C^\infty_c([0,T)\times\Omega_{\delta})$.

\bigskip
\noindent
$\bullet$
Momentum equation
\begin{align} \label{mom_mv}
\int^T_0\int_{\Omega_{\delta}}[\langle Y_{t,x};\mathbf{m}\rangle\partial_t\varphi+\langle Y_{t,x};\frac{\mathbf{m}\otimes\mathbf{m}}{\rho}\rangle:\nabla_x\varphi]dxdt
+\int^T_0\int_{\Omega_{\delta}}\langle Y_{t,x};p(\rho)\rangle\text{div}\varphi dxdt\nonumber\\
=-\int_{\Omega_{\delta}}\langle Y_{0,x};\mathbf{m}\rangle\varphi(0)dx
-\int^T_0\int_{\Omega_{\delta}}\nabla_x\varphi: d\mu_{D}^{M},
\end{align}
for all $\varphi\in C^\infty_c([0,T)\times\Omega_{\delta};\mathbb{R}^3)$ and a signed measure $\mu_{D}^{M}\in\mathcal{M}([0,T]\times\Omega_{\delta};\mathbb{R}^3\times\mathbb{R}^3)$ characterizing the \textit{concentration defect}.

\bigskip
\noindent
$\bullet$
Energy inequality
\begin{align} \label{ei_mv}
\int_{\Omega_{\delta}}\langle &Y_{\tau,x};\frac{1}{2}\frac{|\mathbf{m}|^2}{\rho}+\big{(}P(\rho)-P'(\widetilde{\rho})(\rho-\widetilde{\rho})-P(\widetilde{\rho})\big{)}\rangle dx+\mathcal{D}(\tau)\nonumber\\
&\leq\int_{\Omega_{\delta}}\langle Y_{0,x};\frac{1}{2}\frac{|\mathbf{m}|^2}{\rho}+\big{(}P(\rho)-P'(\widetilde{\rho})(\rho-\widetilde{\rho})-P(\widetilde{\rho})\big{)}\rangle dx
\end{align}
for a.a $\tau\in(0,T)$, where
\begin{align} \label{P}
P(\rho)=\rho\int^\rho_{\widetilde{\rho}}\frac{p(z)}{z^2}dz,
\end{align}
and the non-negative function $\mathcal{D}\in L^\infty(0,T)$ is the so-called \textit{dissipation defect}.

\bigskip
\noindent
$\bullet$
Compatibility conditions
\begin{align} \label{cc}
\int^\tau_0\int_{\Omega_{\delta}}|\mu_{D}^{M}
|dxdt\leq C\int^\tau_0\xi(t)\mathcal{D}(t)dt, \ \ \mbox{for a.a} \ \  \tau\in(0,T), \ \ \mbox{for some} \ \ \xi\in L^1(0,T).
\end{align}

\begin{remark}
The functions
\[
[\rho,\mathbf{m}]\mapsto\frac{\mathbf{m\otimes m}}{\rho}
, \ \ [\rho, \mathbf{m}] \mapsto \frac{ |\mathbf{m}|^2 }{\rho}
\]
are singular on the vacuum set $\rho=0$. We set
\begin{align*}
\frac{|\mathbf{m}|^2}{\rho}=
\left\{
\begin{array}{llll}  \infty,\hspace{5pt}\text{if} \hspace{3pt} \rho=0 \hspace{3pt} and \hspace{3pt} \mathbf{m}\neq0, \\
0,\hspace{5pt}\text{otherwise}.
\end{array}\right.
\end{align*}
Accordingly, it follows from the energy inequality (2.4) that
\begin{align*}
{\rm Supp} [Y_{t,x}]\cap\{[\rho,\mathbf{m}]\in\mathcal{Q}| \ \rho=0,\ \mathbf{m}\neq0]\}=\emptyset \ \mbox{for a.a.}\ (t,x).
\end{align*}
\end{remark}

\begin{remark}
The measure $Y_{0,x}$ plays the role of \emph{initial conditions}.
\end{remark}

\begin{remark}
The proof of an existence of DMV solutions of Euler system was done in the pioneer work  by Neustupa, \cite{Ne}.
Feireisl et al. \cite{j1,j2} proved the existence of \emph{(DMV)} solutions to the non-rotating full Euler system. The existence of DMV solutions to $(1.1)-(1.3)$ can be obtained by analogous methods as in \cite{j2}.
\end{remark}

\subsection{Relative energy inequality}

Motivated by \cite{e6} (see also \cite{e}), we introduce the following modified relative energy functional
\begin{align} \label{ref}
\mathcal{E}(\rho,\mathbf{m}|r, \mathbf{U})=\frac{1}{\delta}
\int_{\Omega_{\delta}}\langle Y_{t,x};\frac{1}{2}\rho|\frac{\mathbf{m}}{\rho}-\mathbf U(t,x)|^2+(P(\rho)-P'(r(t,x))(\rho-r(t,x))-P(r(t,x)))\rangle dx,
\end{align}
where $r>0$, $\mathbf U$ are smooth functions such that $r-\widetilde{\rho}$, $\mathbf U$ are compactly supported in $\Omega_{\delta}$. Any DMV solution of (1.1) satisfies the relative energy inequality (see \cite{e1,e6})
\begin{align} \label{rei}
\mathcal{E}&(\rho,\mathbf{m}|r,\mathbf U)|^{t=\tau}_{t=0}+\mathcal{D(\tau)}\leq
\frac{1}{\delta}\int_{0}^{\tau}\mathcal{R}\left(\varrho,\mathbf{u}\mid{r},\mathbf{U}\right)dt,
\end{align}
with the remainder
\begin{align} \label{rem}
&\mathcal{R}\left(\varrho,\mathbf{u}\mid{r},\mathbf{U}\right)=\int_{\Omega_{\delta}}\langle Y_{t,x};(\partial_t \mathbf U+\frac{\mathbf{m}}{\rho}\nabla_x \mathbf U)(\rho \mathbf U-\mathbf{m})\rangle dx\nonumber\\
&+\frac{1}{\epsilon^{2}}
\int_{\Omega_{\delta}}\langle Y_{t,x};(r-\rho)\partial_tP'(r) -p(\rho)\text{div}\mathbf U
-\mathbf{m}\cdot\nabla_xP'(r)\rangle dx\\
&+\int_{\Omega_{\delta}}\nabla_x \mathbf U:d\mu_{D}^{M}
= \mathcal{R}_1 + \mathcal{R}_2 + \mathcal{R}_3, \nonumber
\end{align}
for a.a. $\tau \in [0,T]$,
and any $r, \mathbf U\in$$C^1([0,T]\times\Omega_{\delta})$, $r-\widetilde{\rho}$, $\mathbf U$ compactly supported in $\Omega_{\delta}$.

\subsection{Solution of the incompressible Euler system}

As shown by Oliver \cite{o}, the incompressible Euler system (1.6)-(1.7) possesses a unique strong solution
\begin{align} \label{Ei}
\mathbf{v}\in C([0,T];W^{m,2}(\mathbb{R}^2)),\hspace{5pt}\Pi\in C([0,T];W^{m,2}(\mathbb{R}^2)),\hspace{8pt}m\geq3,
\end{align}
for any
\begin{align} \label{ic_Ei}
\mathbf{v}_0\in W^{m,2}(\mathbb{R}^2),\hspace{5pt}\text{div}\mathbf{v}_0=0.
\end{align}

\section{Initial-data and acoustic waves}
\label{ILP}

In order to introduce the contribution coming from the acoustic waves, we consider the initial data for the system (\ref{cont}) - (\ref{mom}) in the following form
\begin{equation} \label{id_a}
\overline{\rho_\epsilon (0, \cdot) }= \overline{\rho_{0, \epsilon}} = \widetilde{\rho} + \epsilon s_{0, \epsilon}, \ \
\overline{\mathbf{u}_\epsilon (0, \cdot)} = \overline{\mathbf{u}_{0, \epsilon}}
\end{equation}
for $s_{0}\in L^{\infty}\cap L^{1}\left(\mathbb{R}^{2}\right)$ and $\overline{\mathbf{u}_0 }= \mathbf{v}_0 + \nabla_h \Psi_0,\ {\rm div}_x (\tilde{\rho} \mathbf{v}_0) = 0.$
Here, the acoustic contribution comes from the density perturbation described by $s$, and the gradient of the acoustic potential $\Psi$.
We assume \textit{ill-prepared} initial data, namely
\begin{align} \label{id_conv}
\frac{1}{\delta}\int_{{\Omega}_{\delta}}\langle Y_{0,x}^\epsilon;\frac{1}{2}\rho|\frac{\mathbf{m}}{\rho}-\mathbf{u}_{0, \epsilon} (x)|^2+\frac{1}{\epsilon^2}(P(\rho)-P'(\rho_{0, \epsilon})(\rho- \rho_{0, \epsilon})-P(\rho_{0, \epsilon}))\rangle dx\rightarrow0
\end{align}
as $\epsilon\rightarrow0$.
This could be rephrased as follows
\begin{align} \label{id_conv_1}
\frac{\overline{\rho_{0,\epsilon}}-\widetilde{\rho}}{\epsilon} \ \mbox{bounded in} \ L^{\infty}\left(\mathbb{R}^{2}\right), \ \ \frac{\overline{\rho_{0,\epsilon}}-\widetilde{\rho}}{\epsilon}\rightarrow s_{0} \ \mbox{in} \ L^{1}\left(\mathbb{R}^{2}\right), \ \ \overline{\mathbf{u}_{0,\epsilon}}\rightarrow \overline{\mathbf{u}_0} \ \mbox{in} \ L^{2}\left(\mathbb{R}^{2}, \mathbb{R}^{3}\right).
\end{align}
The corresponding (two-dimensional) acoustic system reads
\begin{equation} \label{as}
\left\{
\begin{array}{llll}  \epsilon\partial_{t}s_\epsilon+\widetilde{\rho}\Delta_h\Psi_\epsilon=0, \\
\epsilon\partial_t\nabla_h\Psi_\epsilon+\frac{p'(\widetilde{\rho})}{\widetilde{\rho}}\nabla_hs_\epsilon=0,
\end{array}\right.
\end{equation}
with the initial data
\begin{equation} \label{id_as}
s_\epsilon (0, \cdot) = s_{0},\ \nabla_h \Psi_\epsilon (0, \cdot) = \nabla_h \Psi_0.
\end{equation}
For technical reasons, the initial data must be smoothed and cut-off via suitable regularization operators, namely
\begin{equation} \label{reg}
s_\epsilon(0,\cdot)=s_{0, \eta}=\frac{\widetilde{\rho}}{p'(\widetilde{\rho})}[\frac{p'(\widetilde{\rho})}{\widetilde{\rho}}s_{0}]_{\eta},\hspace{5pt}
\nabla_h \Psi_\epsilon(0,\cdot)= \nabla_h \Psi_{0, \eta}= \nabla_h [\Psi_{0}]_{\eta},
\end{equation}
where $[\cdot]_\eta$ denotes the regularization.

Denoting the corresponding (smooth) solutions of the acoustic system as  $(s_{\epsilon, \eta}$, $\Psi_{\epsilon, \eta})$, the energy conservation and the standard energy estimates hold in two-dimension, namely
\begin{equation} \label{ec}
\frac{d}{dt}\int_{\mathbb{R}^{2}}\frac{1}{2}\left[p^{\prime}\left(\widetilde{\rho}\right)\left|s_{\epsilon,\eta}(t,\cdot)\right|^{2}+\widetilde{\rho}^{2}\left|\Psi_{\epsilon,\eta}(t,\cdot)\right|^{2}\right]dx=0
\end{equation}
with $p^{\prime}\left(\widetilde{\rho}\right)=a^2$, $a>0$ velocity of sound, and
\begin{align} \label{es}
	\sup_{t\in[0,T]}[\|\Psi_{\epsilon,\eta}(t,\cdot)\|_{W^{m,2}\left(\mathbb{R}^{2}\right)}+\|s_{\epsilon,\eta}(t,\cdot)\|_{W^{m,2}\left(\mathbb{R}^{2}\right)}]\leq C(m,\eta)[\|\nabla_{x}\Psi_{0,\eta}\|_{W^{m,2}\left(\mathbb{R}^{2}\right)}+\|s_{0,\eta}\|_{W^{m,2}\left(\mathbb{R}^{2}\right)}],
\end{align}
for any fixed $m \geq 0$ and $\eta > 0$.

Moreover, the so-called \textit{dispersive estimate} reads as follows (see \cite{sr})
\begin{align} \label{des}
\|\Psi_{\epsilon,\eta}(t,\cdot)\|_{L^{q}\left(\mathbb{R}_{+},W^{k,p}\left(\mathbb{R}^{2}\right)\right)}+\|s_{\epsilon,\eta}(t,\cdot)\|_{L^{q}\left(\mathbb{R}_{+},W^{k,p}\left(\mathbb{R}^{2}\right)\right)}\leq C\epsilon^{\frac{1}{q}}\|\nabla_{x}\Psi_{0,\eta}\|_{W^{m,2}\left(\mathbb{R}^{2}\right)}+\|s_{0,\eta}\|_{W^{m,2}\left(\mathbb{R}^{2}\right)}
\end{align}
for  any
\begin{align} \label{coeff}
&p\in\left(2,\infty\right), \ \ \frac{2}{q}=\frac{1}{2}-\frac{1}{p}, \ \ q\in\left(4,\infty\right)\nonumber\\
&m \geq 0, \ \ k=0,1,...,m-1.
\end{align}

\subsection{Incompressible limit on thin domain -- main result}

Our main result reads as follows.
\begin{theorem} \label{TH}
Let $p\in C^1(0,\infty)\cap C[0,\infty)$ satisfying
$$
p'(\rho)>0 \ \ \mbox{for all} \ \ \rho>0, \ \
\limsup_{\rho\rightarrow{\infty}}
\frac{p(\rho)}{P(\rho)}=P_\infty<\infty,
$$
$$
\liminf_{\rho\rightarrow{\infty}}\frac{p(\rho)}{\rho^\gamma}\geq p_\infty>0 \ \ \mbox{for some} \ \ \gamma>1.
$$
Let $\{Y_{t,x}^\epsilon\}_{(t,x)\in[0,T]\times\Omega}$ be a family of DMV solutions to the scaled compressible Euler system (\ref{cont})-(\ref{mom}) satisfying the compatibility condition (\ref{cc}) with a function $\xi$ independent of $\epsilon$. Let the initial data $\{ Y_{0,x}^\epsilon \}_{x \in \Omega}$ be ill-prepared, namely (\ref{id_conv}) holds.

Then,
\begin{align*}
&\mathcal{D}^\epsilon\rightarrow0\hspace{5pt}\text{in}\hspace{3pt}L^\infty(0,T),\\
&ess\sup_{t\in(\eta,T)}\frac{1}{\delta}\int_0^\delta \int_B\langle Y_{t,x}^\epsilon;\frac{1}{2}\rho|\frac{\mathbf{m}}{\rho}-\mathbf{v}|^2+\frac{1}{\epsilon^2}(P(\rho)-P'(\widetilde{\rho})(\rho- \widetilde{\rho})-P(\widetilde{\rho}))\rangle dx\rightarrow0
\end{align*} as $\epsilon\rightarrow0$ for any compact $B\subset\mathbb{R}^2$ and any $0<\eta<T$, where $\mathbf v$ is the solution of the incompressible Euler system (\ref{cont_E})-(\ref{mom_E}) with initial data $\mathbf v_0=P[\overline{\mathbf u_0}]$, where $P$ denotes the standard Helmholtz projection onto the space of solenoidal functions.
\end{theorem}
The rest of the paper is devoted to the proof of Theorem \ref{TH}.

\section{Energy bounds}

We introduce the following decomposition (see for example \cite{e2})
\begin{align*}
h (\rho, \mathbf{m} ) =[h]_{ess} (\rho, \mathbf{m}) +[h]_{res} (\rho, \mathbf{m})   ,\hspace{10pt}[h]_{ess}=\psi(\rho )h(\rho, \mathbf{m}),\hspace{5pt}[h]_{res}=(1-\psi(\rho))h(\rho, \mathbf{m}),
\end{align*}
where
\begin{align*}
\psi \in C_c^\infty(0,\infty),\hspace{5pt}0\leq\psi(\rho)\leq1,\hspace{5pt}\psi(\rho)=1\hspace{3pt}\text{for all}\hspace{3pt}\rho\in[\frac{1}{2}\min_{\Omega_{\delta}}\widetilde{\rho},2\max_{\Omega_{\delta}}\widetilde{\rho}].
\end{align*}
The right-hand side of the energy inequality (2.4) is bounded uniformly for $\epsilon \to 0$. Consequently, we have
\begin{align} \label{en_b}
ess\sup_{t\in(0,T)}\frac{1}{\delta}\int_{\Omega_{\delta}}\langle Y_{t,x}^\epsilon;\frac{1}{2}\frac{|\mathbf{m}|^2}{\rho}+\frac{1}{\epsilon^2}(P(\rho)-P'(\widetilde{\rho})(\rho- \widetilde{\rho})-P(
\widetilde{\rho}))\rangle dx\leq C.
\end{align}
Now, in order to derive further bounds, we follow the analysis in \cite{e6}. Since
\[
P^{\prime\prime}\left(\rho\right)=\frac{p^{\prime}\left(\varrho\right)}{\varrho}, \ \mbox{for} \ \rho>0,
\]
the convexity of $P\left(\rho\right)$ gives
\begin{align} \label{press_1}
&\left|\rho-\widetilde{\rho}\right|^{2}\leq C\left(\delta\right)\left(P\left(\rho\right)-P^{\prime}\left(\widetilde{\rho}\right)\left(\rho-\widetilde{\rho}\right)-P\left(\widetilde{\rho}\right)\right)
\\
&\mbox{whenever} \ \ 0<\delta\leq\rho, \ \ \widetilde{\rho}\leq\frac{1}{\delta}, \ \ \delta>0; \nonumber
\end{align}
\begin{align} \label{press_2}
&1+\left|\rho-\widetilde{\rho}\right|+P\left(\rho\right)\leq C\left(\delta\right)\left(P\left(\rho\right)-P^{\prime}\left(\widetilde{\rho}\right)\left(\rho-\widetilde{\rho}\right)-P\left(\widetilde{\rho}\right)\right)
\\
&\mbox{if} \ \ 0<2\delta<\widetilde{\rho}<\frac{1}{\delta}, \ \ \varrho\in\left[0,\delta\right)\cup\left[\frac{1}{\delta}, \infty\right), \ \ \delta>0. \nonumber
\end{align}
The above estimates give
\begin{equation} \label{dens_est}
ess\sup_{t\in(0,T)}\frac{1}{\delta}\int_{\Omega_{\delta}}\langle Y_{t,x}^{\epsilon};\left[\left|\frac{\rho-\widetilde{\rho}}{\epsilon}\right|^{2}\right]_{ess}\rangle+\langle Y_{t,x}^{\epsilon};\left[\frac{P(\rho)+1}{\epsilon^{2}}\right]_{res}\rangle dx\leq C.
\end{equation}
Moreover, we have
\begin{equation} \label{mom_est}
\frac{1}{\delta}\int_{\Omega_{\delta}}\langle Y_{t,x}^{\epsilon};\left[\left|\mathbf{m}\right|^{2}\right]_{ess}\rangle dx\leq\frac{1}{\delta}\int_{\Omega_{\delta}}\langle Y_{t,x}^{\epsilon};\left[\rho\left|\frac{\mathbf{m}}{\rho}\right|^{2}\right]_{ess}\rangle dx\leq C.
\end{equation}
From
\[
\left[\left|\mathbf{m}\right|\right]_{res}\leq\frac{\left|\mathbf{m}\right|}{\sqrt{\rho}}\left[\sqrt{\rho}\right]_{res},
\]
we have
\[
\left[\left|\mathbf{m}\right|^{\frac{2\gamma}{\gamma+1}}\right]_{res}\leq C\left(\epsilon\rho\left|\frac{\mathbf{m}}{\rho}\right|^{2}+\frac{1}{\epsilon}\left[\rho^{\gamma}\right]_{res}\right).
\]
Consequently
\begin{equation} \label{mom_est_2}
\frac{1}{\delta}\int_{\Omega_{\delta}}\langle Y_{t,x}^{\epsilon};\left[\left|\mathbf{m}\right|^{\frac{2\gamma}{\gamma+1}}\right]_{res}\rangle dx\leq\epsilon C.
\end{equation}
Finally, we obtain
\begin{align} \label{bounds}
&\langle Y_{t,x}^{\epsilon};\overline{\mathbf{m}}\rangle\hspace{5pt}\text{bounded in}\hspace{3pt}L^{\infty}(0,T;L^{2}(\mathbb{R}^{2})+L^{\frac{2\gamma}{\gamma+1}}(\mathbb{R}^{2})),
\nonumber\\
&\langle Y_{t,x}^{\epsilon};[\frac{\overline{\rho}-\widetilde{\rho}}{\epsilon}]_{ess}\rangle\hspace{5pt}\text{bounded in}\hspace{3pt}L^{\infty}(0,T;L^{2}(\mathbb{R}^{2})),
\nonumber\\
&\epsilon^{-\frac{2}{\gamma}}\langle Y_{t,x}^{\epsilon};[\overline{\rho}]_{res}\rangle\hspace{5pt}\text{bounded in}\hspace{3pt}L^{\infty}(0,T;L^{\gamma}(\mathbb{R}^{2})),
\end{align}
where we averaged in the sense of (\ref{avg}), noticing
\[
\left[\left|\overline{\mathbf{m}}\right|^{2}\right]_{ess}\leq\left[\overline{\left|\mathbf{m}\right|^{2}}\right]_{ess}\hspace{5pt}\text{bounded in}\hspace{3pt}L^{\infty}(0,T;L^{1}(\mathbb{R}^{2}),
\]
\[
\left[\left|\overline{\mathbf{m}}\right|^{\frac{2\gamma}{\gamma+1}}\right]_{res}\leq\left[\overline{\left|\mathbf{m}\right|^{\frac{2\gamma}{\gamma+1}}}\right]_{res}\hspace{5pt}\text{bounded in}\hspace{3pt}L^{\infty}(0,T;L^{1}(\mathbb{R}^{2}),
\]
\[
\left[\left|\frac{\overline{\rho}-\widetilde{\rho}}{\epsilon}\right|^{2}\right]_{ess}
\leq\left[\overline{\left|\frac{\rho-\widetilde{\rho}}{\epsilon}\right|^{2}}\right]_{ess}\hspace{5pt}\text{bounded in}\hspace{3pt}L^{\infty}(0,T;L^{1}(\mathbb{R}^{2}),
\]
\[
ess\sup_{t\in (0,T)} \|\left[\overline{\rho}\right]_{res}\|_{L^{\gamma}(\mathbb{R}^{2})}^{\gamma}\leq ess\sup_{t\in (0,T)} \left\|\left[\overline{\rho^{\gamma}}\right]_{res}\right\|_{L^{1}(\mathbb{R}^{2})}\leq c \epsilon ^{2};
\]
and where we used the following Jensen's inequality
\[
\langle Y_{t,x}^{\epsilon};\left|f\right|\rangle^{p}\leq\langle Y_{t,x}^{\epsilon};\left|f\right|^{p}\rangle, \ \  p\geq1.
\]

\section{Convergence}
The proof of Theorem 3.1 is based on the ansatz $r=\widetilde{\rho}+\epsilon s_{\epsilon,\eta}$ and $\mathbf U=\mathbf V+\nabla_x\Psi_{\epsilon,\eta}$ in the relative energy inequality (2.8). Here, $\mathbf V=(\mathbf v,0)$, with $\mathbf v$ solution of incompressible Euler system (1.6)-(1.7), and $(s_{\epsilon,\eta},\Psi_{\epsilon,\eta})$ solution of the acoustic system (\ref{as}) with $\nabla_x\Psi_{\epsilon,\eta}=(\nabla_h\Psi_{\epsilon,\eta},0)$.
For simplicity, we drop the subscript $\epsilon,\eta$ in $(s_{\epsilon,\eta},\Psi_{\epsilon,\eta})$. Moreover, we assume $\widetilde{\rho}=1$. In the following, we estimate each $\mathcal{R}_j$ ($j=1,2,3$) in (\ref{rem}).

\subsection{The convective term} We write
$$
\frac{1}{\delta}\int_0^\tau\mathcal{R}_1 dt = \frac{1}{\delta}\int_0^\tau\int_{\Omega_{\delta}}\langle Y_{t,x}^{\epsilon}; \left(\partial_{t}\mathbf{U}+\mathbf{U}\cdot\nabla_x\mathbf{U}\right)\cdot\left(\varrho\mathbf{U}-\mathbf{m}\right)\rangle dxdt
$$
\begin{equation}\label{re1_1}
+\frac{1}{\delta}\int_0^\tau\int_{\Omega_\delta} \langle Y_{t,x}^{\epsilon};\left(\frac{\mathbf{m}}{\varrho}-\mathbf{U}\right)\cdot\nabla_x\mathbf{U}\cdot\left(\varrho\mathbf{U}-\mathbf{m}\right)\rangle dxdt.
\end{equation}
The last term is controlled by
\[
\int_0^\tau\|\nabla_h\mathbf{v}(t,\cdot)\|_{L^\infty(\mathbb{R}^{2})}\mathcal{E}(t) dt + \frac{1}{\delta}\int_0^\tau\int_{\Omega_\delta}\langle Y_{t,x}^{\epsilon}; \left(\frac{\mathbf{m}}{\varrho}-\mathbf{U}\right)\cdot\nabla_x{\nabla_x\Psi}\cdot\left(\varrho\mathbf{U}-\mathbf{m}\right) \rangle dxdt
\]
\[
\leq \int_0^\tau c(t)\mathcal{E}(t) dt - \frac{1}{\delta}\int_0^\tau\int_{\Omega_\delta} \langle Y_{t,x}^{\epsilon}; \frac{\mathbf{m}}{\varrho}\otimes\mathbf{m}:\nabla_x{\nabla_x\Psi} \rangle dxdt
\]
\begin{equation}\label{re1_2}
+\frac{2}{\delta}\int_0^\tau\int_{\Omega_\delta} \langle Y_{t,x}^{\epsilon}; \left(\mathbf{m}\otimes\mathbf{U}\right):\nabla_x{\nabla_x\Psi} \rangle dxdt
-  \frac{1}{\delta}\int_0^\tau\int_{\Omega_\delta} \langle Y_{t,x}^{\epsilon}; \varrho\left(\mathbf{U}\otimes\mathbf{U}\right):\nabla_x{\nabla_x\Psi} \rangle dxdt.
\end{equation}
For the $\frac{\mathbf{m}}{\varrho}\otimes\mathbf{m}$ term and from (4.1), we have
\[
\left|\frac{1}{\delta}\int_0^\tau\int_{\Omega_\delta} \langle Y_{t,x}^{\epsilon}; \frac{\mathbf{m}}{\varrho}\otimes\mathbf{m}:\nabla_x{\nabla_x\Psi} \rangle dxdt  \right| \leq c(T)\left\|\overline{\frac{|\mathbf{m}|^2}{\varrho}}\right\|_{L_T^{\infty}(L^1(\mathbb{R}^{2}))}\left\| \nabla_h^2\Psi \right\|_{L_T^8(L^{\infty}(\mathbb{R}^{2}))}
\]
\begin{equation}\label{re1_4}
\leq c(\eta,T)\left\|\overline{\frac{|\mathbf{m}|^2}{\varrho}}\right\|_{L_T^{\infty}(L^1(\mathbb{R}^{2}))}\left\| \nabla_h^2\Psi \right\|_{L_T^8(W^{1,4}(\mathbb{R}^{2}))}\leq c(\eta,T)\epsilon^{\frac{1}{8}}
\end{equation}
Moreover, by using the uniform bound of $\langle Y_{t,x}^{\epsilon};\overline{\mathbf{m}}\rangle$ in $L^{\infty}(0,T;L^2+L^{\frac{2\gamma}{\gamma + 1}}(\mathbb{R}^{2}))$,
\[
\left|\frac{1}{\delta}\int_0^\tau\int_{\Omega_\delta} \langle Y_{t,x}^{\epsilon}; \left(\mathbf{m}\otimes\mathbf{U}\right):\nabla_x{\nabla_x\Psi}\rangle dxdt  \right|
\]
\[
\leq c(T)\left\|\langle Y_{t,x}^{\epsilon};\overline{\mathbf{m}}\rangle\right\|_{L_T^{\infty}(L^2+L^{\frac{2\gamma}{\gamma + 1}}(\mathbb{R}^{2}))} \left\|{\mathbf{U}}\right\|_{L_T^{\infty}(L^4+L^{\frac{6\gamma}{2\gamma-3}}(\mathbb{R}^{3}))} \left\|\nabla_h^2\Psi\right\|_{L_T^{8}(L^4(\mathbb{R}^{2}))+L^6_T(L^{6}(\mathbb{R}^{2}))}
\]
\begin{equation}\label{re1_5}
\leq c(T)c(\eta)\left(\epsilon^{\frac{1}{8}} + \epsilon^{\frac{1}{6}}\right) \leq c(\eta,T)\epsilon^{\frac{1}{8}}.
\end{equation}
For the last $\mathbf{U}\otimes\mathbf{U}$ term in (\ref{re1_2}), we have
\[
\left| \frac{1}{\delta}\int_0^\tau\int_{\Omega_\delta} \langle Y_{t,x}^{\epsilon}; \varrho\left(\mathbf{U}\otimes\mathbf{U}\right):\nabla_x{\nabla_x\Psi} \rangle dxdt  \right|
\]
\[
\leq \epsilon \left| \frac{1}{\delta}\int_0^\tau\int_{\Omega_\delta} \langle Y_{t,x}^{\epsilon}; \frac{\varrho-1}{\epsilon}\left(\mathbf{U}\otimes\mathbf{U}\right):\nabla_x{\nabla_x\Psi} \rangle dxdt \right| +
\left| \int_0^\tau\int_{\mathbb{R}^{2}} \langle Y_{t,x}^{\epsilon}; \left(\mathbf{U}\otimes\mathbf{U}\right):\nabla_x{\nabla_x\Psi} \rangle dxdt \right|
\]
\begin{equation}\label{re1_7}
\leq c(T)\epsilon + c(\eta,T)\epsilon^{\frac{1}{8}} \leq c(\eta,T)\epsilon^{\frac{1}{8}}.
\end{equation}
For the first term on the right side of (\ref{re1_1}),
$$
\frac{1}{\delta}\int_0^\tau\int_{\Omega_\delta}\langle Y_{t,x}^{\epsilon}
; \left(\partial_{t}\mathbf{U}+\mathbf{U}\cdot\nabla_{x}\mathbf{U}\right)\cdot\left(\varrho\mathbf{U}-\mathbf{m}\right) \rangle dxdt
$$
$$
=\frac{1}{\delta}\int_0^\tau\int_{\Omega_\delta}\langle Y_{t,x}^{\epsilon}
; \left(\partial_{t}\mathbf{V}+\mathbf{V}\cdot\nabla_x\mathbf{V}\right)\cdot\left(\varrho\mathbf{U}-\mathbf{m}\right) \rangle dxdt
+ \frac{1}{\delta}\int_0^\tau\int_{\Omega_\delta}\langle Y_{t,x}^{\epsilon}
; \partial_{t}\nabla_x{\Psi}\cdot\left(\varrho\mathbf{U}-\mathbf{m}\right) \rangle dxdt
$$
$$
+\frac{1}{\delta}\int_0^\tau\int_{\Omega_\delta} \langle Y_{t,x}^{\epsilon}
; \nabla_x{\Psi}\cdot\nabla_x\nabla_x{\Psi}\cdot\left(\varrho\mathbf{U}-\mathbf{m}\right) \rangle dxdt
$$
\begin{equation}\label{re3}
+\frac{1}{\delta}\int_0^\tau\int_{\Omega_\delta} \langle Y_{t,x}^{\epsilon}
; \left(\mathbf{V}\cdot\nabla_x(\nabla_x\Psi) + \nabla_x\Psi\cdot\nabla_x\mathbf{V}\right)\cdot\left(\varrho\mathbf{U}-\mathbf{m}\right) \rangle dxdt.
\end{equation}
Since $\mathbf V=(\mathbf v,0)$ and $\mathbf{v}$ is the solution to the Euler equations (\ref{mom_E}), we have
\[
\frac{1}{\delta}\int_0^\tau\int_{\Omega_\delta}\langle Y_{t,x}^{\epsilon}
; \left(\partial_{t}\mathbf{V}+\mathbf{V}\cdot\nabla_x\mathbf{V}\right)\cdot\left(\varrho\mathbf{U}-\mathbf{m}\right) \rangle dxdt=
{I_{1}+I_{2}},
\]
where
\[
{I_{1}}=\frac{1}{\delta}\int_0^\tau\int_{\Omega_\delta}  \langle Y_{t,x}^{\epsilon}
; \mathbf{m} \cdot \nabla_x\Pi \rangle dxdt =\frac{1}{\delta}\int_{\Omega_\delta}  \langle Y_{t,x}^{\epsilon}
; \varrho  \Pi dx\rangle dx|_{t=0}^\tau
- \frac{1}{\delta}\int_0^\tau\int_{\Omega_\delta}  \langle Y_{t,x}^{\epsilon}
; \varrho \partial_t \Pi \rangle dxdt
\]
\begin{equation}\label{d3}
=\epsilon \frac{1}{\delta}\int_{\Omega_\delta} \langle Y_{t,x}^{\epsilon}
; \frac{\varrho - 1}{\epsilon} \Pi dx\left|_{t=0}^\tau\right \rangle dx
- \epsilon \frac{1}{\delta}\int_0^\tau\int_{\Omega_\delta}  \langle Y_{t,x}^{\epsilon}
; \frac{\varrho -1}{\epsilon} \partial_t \Pi \rangle dxdt \leq c(\eta,T)\epsilon,
\end{equation}
with $\nabla_x\Pi=(\nabla_h\Pi,0)$ and $\Pi$ the pressure in $(1.7)$, and
	\[
	|{I_{2}}|=\left|\frac{1}{\delta}\int_{0}^{\tau}\int_{\Omega_{\delta}}
	\langle Y_{t,x}^{\epsilon}
	; \varrho\mathbf{U}\cdot\nabla_x\Pi \rangle dxdt\right|\leq\left|\frac{1}{\delta}\int_{0}^{\tau}\int_{\Omega_{\delta}}
	\langle Y_{t,x}^{\epsilon}
	; \left(\varrho-1\right)\cdot\mathbf{U}\cdot\nabla_x\Pi \rangle dxdt \right|
	\]
	\begin{equation}
	+\left|\frac{1}{\delta}\int_{0}^{\tau}\int_{\Omega_{\delta}}
	\langle Y_{t,x}^{\epsilon}
	; \mathbf{U}\cdot\nabla_x\Pi \rangle dxdt \right|.\label{split}
	\end{equation}
	Similarly to the analysis above, for the first term on the right hand side of (\ref{split}), we have
	\[
	\left|\frac{1}{\delta}\int_{0}^{\tau}\int_{\Omega_{\delta}} \langle Y_{t,x}^{\epsilon}
	; \left(\varrho-1\right)\cdot\mathbf{U}\cdot\nabla_x\Pi \rangle dxdt \right|
	\leq \epsilon\left|\frac{1}{\delta}\int_{0}^{\tau}\int_{\Omega_{\delta}} \langle Y_{t,x}^{\epsilon}
	; \frac{\left(\varrho-1\right)}{\epsilon}\cdot\mathbf{U}\cdot\nabla_x\Pi \rangle dxdt \right|
	\]
	\[
	\leq c(T)\epsilon
	\]
	For the second term on the right hand side of (\ref{split}), we have
	\begin{equation}
	\frac{1}{\delta}\int_{0}^{\tau}\int_{\Omega_{\delta}} \langle Y_{t,x}^{\epsilon}
	; \mathbf{U}\cdot\nabla_x\Pi  \rangle dxdt=\frac{1}{\delta}\int_{0}^{\tau}\int_{\Omega_{\delta}} \langle Y_{t,x}^{\epsilon}
	; \mathbf{V}\cdot\nabla_x\Pi \rangle dxdt+\frac{1}{\delta}\int_{0}^{\tau}\int_{\Omega_{\delta}} \langle Y_{t,x}^{\epsilon}
	; \nabla_x\Psi\cdot\nabla_x\Pi \rangle dxdt.\label{press_conv2-4}
	\end{equation}
	Performing integration by parts in the first term on the right-hand
	side of (\ref{press_conv2-4}), we have
	\[
	\frac{1}{\delta}\int_{0}^{\tau}\int_{\Omega_{\delta}} \langle Y_{t,x}^{\epsilon}
	; \textrm{div}_x\mathbf{V}\Pi \rangle dxdt=0
	\]
	thanks to incompressibility condition, $\textrm{div}_x\mathbf{V}=0$.
	For the second term on the right-hand side of (\ref{press_conv2-4})
	using integration by parts and acoustic equation, we
	have
	\[
	\frac{1}{\delta}\int_{0}^{\tau}\int_{\Omega_{\delta}} \langle Y_{t,x}^{\epsilon}
	; \nabla_x\Psi\cdot\nabla_x\Pi \rangle dxdt=-\frac{1}{\delta}\int_{0}^{\tau}\int_{\Omega_{\delta}} \langle Y_{t,x}^{\epsilon};\Delta_x\Psi\Pi \rangle dxdt
	\]
	\[
	=\epsilon\frac{1}{\delta}\int_{0}^{\tau}\int_{\Omega_{\delta}} \langle Y_{t,x}^{\epsilon}; \partial_{t}s\Pi \rangle dxdt
	\]
	\begin{equation}
	=\epsilon\left[\frac{1}{\delta}\int_{\Omega_{\delta}} \langle Y_{t,x}^{\epsilon}; s\Pi \rangle dx\right]_{t=0}^{t=\tau}-\epsilon\frac{1}{\delta}\int_{0}^{\tau}\int_{\Omega_{\delta}} \langle Y_{t,x}^{\epsilon}; s\partial_{t}\Pi \rangle dxdt,\label{phi_p}
	\end{equation}
that it goes to zero for $\epsilon\rightarrow0$. Moreover, by using similar argument as above, the last two terms in (\ref{re3}) are of order
\begin{equation}\label{re5}
c(\eta,T)(1+\epsilon) \|\nabla_h\Psi \|_{L^8_T(W^{1,4}(\mathbb{R}^{2}))} \leq c(\eta,T)\epsilon^{\frac{1}{8}}.
\end{equation}
Finally, using ${\rm div}_x\mathbf{V}=0$, we get
$$
\frac{1}{\delta}\int_0^\tau\int_{\Omega_\delta} \langle Y_{t,x}^{\epsilon}; \partial_{t}\nabla_x{\Psi}\cdot\left(\varrho \mathbf{U}-\mathbf{m}\right) \rangle dxdt = - \frac{1}{\delta}\int_0^\tau\int_{\Omega_\delta} \langle Y_{t,x}^{\epsilon}; \mathbf{m}\cdot \partial_{t}\nabla_x{\Psi} \rangle dxdt
$$
\begin{equation}\label{re6}
+ \frac{1}{\delta}\int_0^\tau\int_{\Omega_\delta} \langle Y_{t,x}^{\epsilon}; (\varrho - 1) \mathbf{V}\cdot\partial_{t}\nabla_x\Psi \rangle dxdt + \frac{1}{\delta}\int_0^\tau\int_{\Omega_\delta} \langle Y_{t,x}^{\epsilon}; \varrho \partial_{t} \nabla_x\Psi \cdot \nabla_x\Psi \rangle dxdt
\end{equation}
The first term on the right side of (\ref{re6}) will be cancelled later by the pressure term while, by using the acoustic wave equations (\ref{as}),
the second term equals to
$$
\frac{1}{\delta}\int_0^\tau\int_{\Omega_\delta} \langle Y_{t,x}^{\epsilon}; \frac{\varrho-1}{\epsilon} \epsilon \partial_{t} \nabla_x{\Psi}\cdot \mathbf{V} \rangle dxdt = -\frac{1}{\delta}\int_0^\tau\int_{\Omega_\delta} \langle Y_{t,x}^{\epsilon}; \frac{\varrho-1}{\epsilon} a^2\nabla_x s\cdot \mathbf{V} \rangle dxdt
$$
\[
\leq c(T)\left\|\frac{\overline{\varrho}-1}{\epsilon}\right\|_{L^\infty_T(L^2+L^{\gamma_2}(\mathbb{R}^{2}))} \left\|\mathbf{v}\right\|_{L^\infty_T(L^4+L^{\frac{4\gamma}{3\gamma-4}}(\mathbb{R}^{2}))} \left\|\nabla_h s\right\|_{L^{8}_T(L^4+L^{4}(\mathbb{R}^{2}))}
\]
\begin{equation}\label{re7}
\leq c(\eta,T)\epsilon^{\frac{1}{8}}, \, \gamma_2= \min\{2,\gamma\}.
\end{equation}
Finally, by using the acoustic equations, $\epsilon\partial_t\nabla_x\Psi = -a^{2} \nabla_x s$,
\[
\frac{1}{\delta}\int_0^\tau\int_{\Omega_\delta} \langle Y_{t,x}^{\epsilon}; \varrho \partial_{t} \nabla_x\Psi \cdot \nabla_x\Psi \rangle dxdt
\]
\[
= - a^{2}\frac{1}{\delta}\int_0^\tau\int_{\Omega_\delta} \langle Y_{t,x}^{\epsilon}; \frac{\varrho-1}{\epsilon}  \nabla_x s \cdot \nabla_x \Psi \rangle dxdt + \frac{1}{2}\int_{\mathbb{R}^{2}} \langle Y_{t,x}^{\epsilon}; |\nabla_h \Psi|^2 \left|_{t=0}^{\tau}\right. \rangle dx
\]
\begin{equation}\label{re9}
\leq c(\eta,T)\epsilon^{\frac{1}{8}} + \frac{1}{2}\int_{\mathbb{R}^{2}} \langle Y_{t,x}^{\epsilon}; |\nabla_h \Psi|^2 dx\left|_{t=0}^{\tau}\right. \rangle dx.
\end{equation}
From (\ref{re1_1}) to (\ref{re9}) we find
\begin{equation}\label{r1final}
\frac{1}{\delta}\int_0^\tau \mathcal{R}_1 dt \leq c(\eta,T)\epsilon^{\frac{1}{8}} + \int_0^\tau c(t)\mathcal{E}(t) dt + \frac{1}{2}\int_{\mathbb{R}^{2}} \langle Y_{t,x}^{\epsilon}; |\nabla_h \Psi|^2 \left|_{t=0}^{\tau}\right. \rangle dx - \frac{1}{\delta}\int_0^\tau\int_{\Omega_\delta} \langle Y_{t,x}^{\epsilon}; \varrho\mathbf{u}\cdot \partial_{t}\nabla_x {\Psi} \rangle dxdt.
\end{equation}

\subsection{Terms depending on the pressure}
We recall that
\[
\frac{1}{\delta}\int_0^\tau\mathcal{R}_2 dt = \frac{1}{\epsilon^{2}}\frac{1}{\delta}\int_0^\tau\int_{\Omega_{\delta}} \langle Y_{t,x}^{\epsilon}; \left(\varrho-r\right)\partial_{t}P^{\prime}\left(r\right) - p\left(\varrho \right)\textrm{div}_x\mathbf{U} - \mathbf{m}\cdot \nabla_x P^{\prime}\left(r\right) \rangle dxdt,
\]
where $r = 1 + \epsilon s$. We have
\[
\frac{1}{\epsilon^{2}}\frac{1}{\delta}\int_0^\tau\int_{\Omega_{\delta}} \langle Y_{t,x}^{\epsilon};  \mathbf{m}\cdot \nabla_x P^{\prime}\left(r\right) \rangle dxdt =\frac{1}{\epsilon}\frac{1}{\delta}\int_0^\tau\int_{\Omega_{\delta}}  \langle Y_{t,x}^{\epsilon}; \mathbf{m}\cdot \nabla_x s  {P}^{\prime\prime}(r) \rangle dxdt
\]
\[
=\frac{1}{\delta}\int_0^\tau\int_{\Omega_{\delta}} \langle Y_{t,x}^{\epsilon};  \mathbf{m}\cdot \nabla_x s \frac{ {P}^{\prime\prime}(1+\epsilon s) - {P}^{\prime\prime}(1)}{\epsilon} \rangle dxdt + \frac{1}{\delta}\frac{1}{\epsilon}\int_0^\tau\int_{\Omega_{\delta}} \langle Y_{t,x}^{\epsilon}; a^2 \mathbf{m}\cdot \nabla_x s \rangle dxdt
\]
since ${P}^{\prime\prime}(1)=p'(1)=a^2$. Realizing that
$$\left| \frac{ {P}^{\prime\prime}(1+\epsilon s) - {P}^{\prime\prime}(1)}{\epsilon}\right|\leq c | s |, $$
the first term on the right side is controlled by
\[
c(\eta,T) \left\| \langle Y_{t,x}^{\epsilon};\overline{\mathbf{m}}\rangle\right\|_{L^\infty_T(L^2+L^{\frac{2\gamma}{\gamma+1}}(\mathbb{R}^{2}))} \left\|s \right\|_{L^\infty_T(L^{4}+L^\infty(\mathbb{R}^{2}))} \left\|\nabla_h s\right\|_{L^{8}_T(L^{4}(\mathbb{R}^{2}))+L^{4\gamma}_T(L^{\frac{2\gamma}{\gamma -1}}(\mathbb{R}^{2}))}
\]
\begin{equation}\label{pr1}
\leq c(\eta,T) \epsilon^{\min\{ \frac{1}{8}, \frac{1}{4\gamma} \}}.
\end{equation}
By using the acoustic equations,
\[
\frac{1}{\delta}\frac{1}{\epsilon}\int_0^\tau\int_{\Omega_{\delta}} \langle Y_{t,x}^{\epsilon}; a^2 \mathbf{m}\cdot \nabla_x s  \rangle dxdt
= - \frac{1}{\delta}\int_0^\tau\int_{\Omega_{\delta}}  \langle Y_{t,x}^{\epsilon}; \mathbf{m}\cdot \partial_t\Psi \rangle  dxdt,
\]
which cancels the same term appeared on the right side of (\ref{re6}). Now we write
\[
\frac{1}{\varepsilon^{2}}\frac{1}{\delta}\int_0^\tau\int_{\Omega_{\delta}} \langle Y_{t,x}^{\epsilon}; \left(\varrho-r\right)\partial_{t}P^{\prime}\left(r\right) - p\left(\varrho \right)\textrm{div}_x\mathbf{U} \rangle dxdt
\]
\[
= \frac{1}{\delta}\int_0^\tau\int_{\Omega_{\delta}} \langle Y_{t,x}^{\epsilon}; \frac{\varrho- 1}{\epsilon} P^{\prime\prime}(r)\partial_t s \rangle dxdt + \int_0^\tau\int_{\mathbb{R}^{2}} \langle Y_{t,x}^{\epsilon}; s P^{\prime\prime}(r)\partial_t s \rangle dx_h dt
\]
\[
- \frac{1}{\delta}\int_0^\tau\int_{\Omega_{\delta}} \langle Y_{t,x}^{\epsilon}; \frac{ p (\varrho)-p'(1)(\varrho-1)-p(1)}{\epsilon^2} \Delta_x {\Psi} \rangle dxdt
\]
\begin{equation}\label{pr3}
- \frac{1}{\delta}\int_0^\tau\int_{\Omega_{\delta}} \langle Y_{t,x}^{\epsilon};  p'(1)\frac{\varrho -1}{\epsilon}\frac{1}{\epsilon}\Delta_x \Psi \rangle dxdt.
\end{equation}
Note that
\[
\frac{1}{\delta}\int_0^\tau\int_{\Omega_{\delta}} \langle Y_{t,x}^{\epsilon}; \frac{\varrho- 1}{\epsilon} P^{\prime\prime}(r)\partial_t s \rangle dxdt =\frac{1}{\delta}\int_0^\tau\int_{\Omega_{\delta}} \langle Y_{t,x}^{\epsilon}; \frac{\varrho- 1}{\epsilon} P^{\prime\prime}(1)\partial_t s \rangle dxdt
\]
\[
+ \frac{1}{\delta}\int_0^\tau\int_{\Omega_{\delta}} \langle Y_{t,x}^{\epsilon}; \frac{\varrho- 1}{\epsilon} \left(P^{\prime\prime}(r)-P^{\prime\prime}(1)\right)\partial_t s \rangle dxdt .
\]
We find the first term on the right side is cancelled by the last term in (\ref{pr3}), while the remaining term equals to
\[
- \frac{1}{\delta}\int_0^\tau\int_{\Omega_{\delta}} \langle Y_{t,x}^{\epsilon}; \frac{\varrho- 1}{\epsilon} \frac{P^{\prime\prime}(r)-P^{\prime\prime}(1)}{\epsilon}\Delta_x \Psi \rangle dxdt
\]
\[
\leq c(T)\left\|\frac{\overline{\varrho} -1}{\epsilon}\right\|_{L^\infty_T(L^{2}+L^{\gamma}(\mathbb{R}^{2}))} \left\| s \right\|_{L^\infty_T(L^{4}+L^{\frac{4\gamma}{3\gamma -4}}(\mathbb{R}^{2}))} \left\|\Delta_h \Psi\right\|_{L^{8}_T(L^4+L^4(\mathbb{R}^{2}))}
\]
\begin{equation}\label{pr4}
\leq c(\eta,T)\epsilon^{\frac{1}{8}}.
\end{equation}
Similarly,
\[
\int_0^\tau\int_{\mathbb{R}^{2}}  \langle Y_{t,x}^{\epsilon}; s P^{\prime\prime}(r)\partial_t s \rangle dxdt =  \int_0^\tau\int_{\mathbb{R}^{2}} \langle Y_{t,x}^{\epsilon}; s P^{\prime\prime}(1) \partial_t s \rangle dx_h dt
\]
\[
+ \int_0^\tau\int_{\mathbb{R}^{2}}  \langle Y_{t,x}^{\epsilon}; s \left( P^{\prime\prime}(r)-P^{\prime\prime}(1)\right)\partial_t s \rangle dx_h dt
\]
\[
\leq \frac{1}{2} \int_{\mathbb{R}^{2}} \langle Y_{t,x}^{\epsilon}; a^2 \left| s \right|^2  \left|_{t=0}^\tau \rangle \right. dx + c(T)\left\|s \right\|_{L^\infty_T(L^2(\mathbb{R}^{2}))} \left\|s \right\|_{L^\infty_T(L^4(\mathbb{R}^{2}))} \left\|\Delta_h \Psi\right\|_{L^{8}_T(L^4(\mathbb{R}^{2}))}
\]
\begin{equation}\label{pr6}
\leq \frac{1}{2} \int_{\mathbb{R}^{2}} \langle Y_{t,x}^{\epsilon}; a^2 \left| s \right|^2  \left|_{t=0}^\tau \right. \rangle dx_h + c(\eta,T)\epsilon^{\frac{1}{8}}.
\end{equation}
Finally, realizing that $\frac{1}{\delta}\frac{ p (\varrho)-p'(1)(\varrho-1)-p(1)}{\epsilon^2}$ is uniformly bounded in $L^\infty(0,T;L^1(\Omega_{\delta}))$,
\[
\frac{1}{\delta}\int_0^\tau\int_{\Omega_{\delta}} \langle Y_{t,x}^{\epsilon}; \frac{ p (\varrho)-p'(1)(\varrho-1)-p(1)}{\epsilon^2} \Delta_x {\Psi} \rangle dxdt
\]
\begin{equation}\label{pr8}
\leq c(T)\left\|\Delta_h\Psi\right\|_{L^8_T(L^\infty(\mathbb{R}^{2}))} \leq c(T)\left\|\nabla_h\Psi\right\|_{L^8_T(W^{2,4}(\mathbb{R}^{2}))} \leq c(\eta,T)\epsilon^{\frac{1}{8}}.
\end{equation}
From (\ref{pr1}) to (\ref{pr8}) we conclude that
\begin{equation}\label{prfinal}
\frac{1}{\delta}\int_0^\tau\mathcal{R}_2 dt \leq \frac{1}{2} \int_{\mathbb{R}^{2}} \langle Y_{t,x}^{\epsilon}; a^2 \left| s \right|^2 \left|_{t=0}^\tau \right. \rangle dx_h + c(\eta,T)\epsilon^{\alpha}, \, \alpha=\min\{\frac{1}{8},\frac{1}{4\gamma}\}.
\end{equation}

\subsection{Concentration measure term}

Finally, we can control the concentration measure through the use of the compatibility conditions (\ref{cc}), namely
\[
\frac{1}{\delta}\int_0^\tau\mathcal{R}_3 dt=\frac{1}{\delta}\int_{0}^{\tau}\int_{\Omega_{\delta}}\nabla_{x}\mathbf{U}:\mu_{D}^{M}dxdt\leq\left\Vert \nabla_{x}\mathbf{U}\right\Vert _{L^{\infty}\left(\mathbb{R}^{3}\right)}\int_{0}^{\tau}\xi\left(t\right)\mathcal{D}\left(t\right)dt.
\]

Using the conservation of energy for acoustic wave system and all estimates in the above three subsections, we find
\[
\mathcal{E}\left(\varrho,\mathbf{u}\mid {r},\mathbf{U}\right)\left(\tau\right) + \mathcal{D}(\tau)\leq c(\eta,T)\epsilon^{\alpha} + \int_0^\tau c(t) \mathcal{E}(t) dt + \int_0^\tau \xi(t) \mathcal{D}(t) dt,
\]
where $c(t)=\|\nabla_h\mathbf{v}(t,\cdot)\|_{L^\infty(\mathbb{R}^{2})}\leq c\|\mathbf{v}(t,\cdot)\|_{W^{3,2}(\mathbb{R}^{2})}$ according to Sobolev's embedding lemma. By Gronwall's inequality,
\begin{equation}\label{cpt1}
\mathcal{E}\left(\varrho,\mathbf{u}\mid {r},\mathbf{U}\right)\left(\tau\right) \leq c(\eta,T)\epsilon^\alpha + c(T) \mathcal{E}(0), \, \text{ a.e. } \tau\in (0,T),
\end{equation}
where $c(T) = \exp{\int_0^T\|\nabla_h\mathbf{v}(t,\cdot)\|_{L^\infty(\mathbb{R}^{2})}dt}$. Sending $\epsilon\to 0$ and then $\eta\to 0$ we find
\[
\lim_{\eta\to 0}\lim_{\epsilon\to 0}\mathcal{E}\left(\varrho_{\epsilon},\mathbf{u}_{\epsilon}\mid {{r}_{\epsilon,\eta}},\mathbf{U_{\epsilon,\eta}}\right)\left(\tau\right)=0 \text{ uniformly in }\tau\in (0,T),
\]
as well as
\[
\lim_{\epsilon\to 0}\mathcal{E}\left(\varrho_{\epsilon},\mathbf{u}_{\epsilon}\mid {r}_{\epsilon},\mathbf{U}_{\epsilon}\right)\left(\tau\right) = 0\text{ uniformly in }\tau\in (0,T),
\]
where $r_{\epsilon} = 1+ s_{\epsilon},\,\mathbf{U}_{\epsilon}=(\mathbf{v} + \nabla_h\Psi_{\epsilon},0)$. We thus conclude the proof of Theorem 3.1 by realizing that $\nabla_h\Psi_{\epsilon,\eta}\to 0$ in $L^q(0,T;L^p(\mathbb{R}^{2}))$ as $\epsilon\to 0$ for any $p>2, \ q>4$ according to the dispersive estimates (\ref{des}). Indeed, for any compact set $K\subset\mathbb{R}^{2}$,
\[
\left\|\overline{\sqrt{\varrho_{\epsilon}}\mathbf{u}_{\epsilon}} - \mathbf{v}\right\|_{L^2_T(L^2(K))} \leq \left\|\overline{\sqrt{\varrho_{\epsilon}}\mathbf{u}_{\epsilon}}
- \mathbf{U}_{\epsilon,\eta}\right\|_{L^2_T(L^2(\mathbb{R}^{2}))} + c(T,K)\left\|{\nabla_h{\Psi}_{\epsilon,\eta}} \right\|_{L^q_T(L^p(K))},
\]
which vanishes as $\epsilon\to 0$ and then $\eta\to 0$.

\vskip 0.5cm
\noindent {\bf Acknowledgements}

\vskip 0.1cm

The research of \v S.N. leading to these results has received funding from the Czech Sciences Foundation (GA\v CR), GA19-04243S. The research of T.T. is supported by the NSFC Grant No. 11801138. The research of B.D. is partially supported by the ANR project INFAMIE (ANR-15-CE40-0011).


\end{document}